\documentclass[11pt]{amsart}
\usepackage{amsmath,amssymb,mathrsfs}
\newtheorem{theorem}{Theorem}
\newtheorem{corollary}[theorem]{Corollary}
\newtheorem{lemma}[theorem]{Lemma}
\newtheorem{proposition}[theorem]{Proposition}
\newtheorem*{conjecture}{Conjecture}

\begin{document}

\title{Uniqueness of gradient Ricci solitons}
\author{Simon Brendle}
\address{Department of Mathematics \\ Stanford University \\ 450 Serra Mall, Bldg 380 \\ Stanford, CA 94305} 
\thanks{The author was supported in part by the National Science Foundation under grant DMS-0905628.}
%\begin{abstract}
%We show that a three-dimensional steady gradient Ricci soliton which is asymptotic to the Bryant soliton in a suitable sense must be isometric to the Bryant soliton.
%\end{abstract}
\maketitle

\section{Introduction}

The Ricci flow, introduced by R.~Hamilton \cite{Hamilton} in 1982, has been studied intensively in recent years. In particular, the Ricci flow plays a key role in Perelman's proof of the Poincar\'e conjecture (cf. \cite{Perelman1}, \cite{Perelman2}, \cite{Perelman3}). The Ricci flow also features prominently in the proof of the Differentiable Sphere Theorem for pointwise $1/4$-pinched manifolds (cf. \cite{Brendle}, \cite{Brendle-Schoen1}, \cite{Brendle-Schoen2}, \cite{Brendle-Schoen3}). For an introduction to Ricci flow, see e.g. \cite{Brendle-book} or \cite{Topping}.

In this paper, we are interested in self-similar solutions to the Ricci flow. Such solutions are referred to as Ricci solitons, and were first studied by Hamilton \cite{Hamilton-survey}. Recall that a Riemannian manifold $(M,g)$ is called a steady Ricci soliton if 
\[\text{\rm Ric} + \frac{1}{2} \, \mathscr{L}_\xi(g) = 0\] 
for some vector field $\xi$. Moreover, if $\xi = -\nabla f$ for some smooth function $f: M \to \mathbb{R}$, then $(M,g)$ is referred to as a steady gradient Ricci soliton. Ricci solitons play a fundamental role in the formation of singularities, and have been studied by many authors; see \cite{Cao-survey} for a survey.

The simplest example of a steady Ricci soliton is the so-called cigar soliton in dimension $2$. The cigar soliton is rotationally symmetric, has positive Gaussian curvature, and is asymptotic to a cylinder near infinity. R.~Bryant \cite{Bryant} has constructed an example of a steady gradient Ricci soliton in dimension $3$. This solution is rotationally symmetric and has positive sectional curvature. Bryant's construction can be adapted to higher dimensions. In fact, for each $n \geq 3$, there exists an $n$-dimensional steady gradient Ricci soliton, which is rotationally symmetric and has positive curvature operator. This will be referred to as the Bryant soliton. Other examples of steady Ricci solitons were constructed by H.D.~Cao \cite{Cao} and T.~Ivey \cite{Ivey}.

It was shown by Hamilton that any two-dimensional gradient soliton is isometric to the cigar soliton up to scaling. In \cite{Perelman1}, G.~Perelman conjectured a similar uniqueness property in dimension $3$:

\begin{conjecture}[G.~Perelman \cite{Perelman1}]
Any three-dimensional steady gradient Ricci soliton with positive sectional curvature which satisfies a non-collapsing assumption at infinity is isometric to the Bryant soliton up to scaling.
\end{conjecture} 

We note that H.~Guo \cite{Guo} has obtained interesting results on the asymptotic geometry of a Ricci soliton near infinity. In a recent paper \cite{Cao-Chen}, H.D.~Cao and Q.~Chen proved uniqueness under the additional assumption that $(M,g)$ is locally conformally flat. The same result was proved independently by Catino and Mantegazza \cite{Catino-Mantegazza} under the assumption that $n \geq 4$. 

\begin{theorem}[H.D.~Cao, Q.~Chen \cite{Cao-Chen}]
Let $(M,g)$ be a steady gradient Ricci soliton of dimension $n \geq 3$. If $(M,g)$ is locally conformally flat, then $(M,g)$ is either flat or rotationally symmetric.
\end{theorem}

Throughout this paper, we will assume that $(M,g)$ is a three-dimensional steady gradient Ricci soliton. We will show that $(M,g)$ is rotationally symmetric, provided that $(M,g)$ satisfies certain asymptotic conditions near infinity. To that end, we fix a smooth function $\psi: (0,1) \to \mathbb{R}$ so that $\nabla R + \psi(R) \, \nabla f = 0$ on the Bryant soliton. Moreover, we define 
\[u(s) = \log \psi(s) + \int_{\frac{1}{2}}^s \Big ( \frac{3}{2(1-t)} - \frac{1}{(1-t) \psi(t)} \Big ) \, dt.\] 
Then we have the following result:

\begin{theorem}
\label{main.theorem}
Let $(M,g)$ be a three-dimensional steady Ricci soliton. Suppose that the scalar curvature of $(M,g)$ is positive and approaches zero at infinity. Moreover, we assume that there exists an exhaustion of $M$ by bounded domains $\Omega_l$ such that 
\begin{equation} 
\label{assumption}
\lim_{l \to \infty} \int_{\partial \Omega_l} e^{u(R)} \, \langle \nabla R + \psi(R) \, \nabla f,\nu \rangle = 0. 
\end{equation}
Then $(M,g)$ is rotationally symmetric.
\end{theorem}

The proof of Theorem \ref{main.theorem} is inspired in part by D.C.~Robinson's proof of the uniqueness of the Schwarzschild black hole (cf. \cite{Israel} and \cite{Robinson}).

The author would like to thank Fernando Marques for discussions. He is grateful to the referee for many useful comments.

\section{The key identities}

Let $(M,g)$ be a three-dimensional steady gradient Ricci soliton, so that $\text{\rm Ric} = D^2 f$ for some real-valued function $f$. We first collect some well-known facts:

\begin{proposition}
We have 
\begin{equation} 
\label{gradient.of.R}
\partial_i R = -2 \, \text{\rm Ric}_{ij} \, \partial_j f. 
\end{equation} 
and 
\begin{equation} 
\label{Laplacian.of.R}
\Delta R + 2 \, |\text{\rm Ric}|^2 = -\langle \nabla f,\nabla R \rangle. 
\end{equation}
\end{proposition} 

\textbf{Proof.} 
Using the contracted second Bianchi identity, we obtain 
\begin{align*} 
0 &= \partial_i R - 2 \, g^{kl} \, D_i \text{\rm Ric}_{kl} + 2 \, g^{kl} \, D_k \text{\rm Ric}_{il} \\ 
&= \partial_i R - 2 \, g^{kl} \, D_{i,k,l}^3 f + 2 \, g^{kl} \, D_{k,i,l}^3 f \\ 
&= \partial_i R + 2 \, g^{kl} \, R_{ikjl} \, \partial^j f \\ 
&= \partial_i R + 2 \, \text{\rm Ric}_{ij} \, \partial^j f. 
\end{align*} 
This proves (\ref{gradient.of.R}). To prove (\ref{Laplacian.of.R}), we take the divergence on both sides of the previous identity. This yields 
\begin{align*} 
0 &= \Delta R + 2 \, D^i \text{\rm Ric}_{ij} \, \partial^j f + 2 \, |\text{\rm Ric}|^2 \\ 
&= \Delta R + \partial_j R \, \partial^j f + 2 \, |\text{\rm Ric}|^2, 
\end{align*} 
as claimed. \\

It follows from (\ref{gradient.of.R}) that the sum $R + |\nabla f|^2$ is constant. By scaling, we may ssume that $R + |\nabla f|^2 = 1$ at each point in $M$. We next define a tensor $B_{ijk}$ by 
\begin{align*} 
B_{ijk} 
&= \text{\rm Ric}_{ik} \, \partial_j f - \text{\rm Ric}_{ij} \, \partial_k f \\ 
&- \frac{1}{4} \, \big ( (\partial_j R + 2R \, \partial_j f) \, g_{ik} - (\partial_k R + 2R \, \partial_k f) \, g_{ij} \big ). 
\end{align*} 
Note that the tensor $B_{ijk}$ vanishes on the set $\{R=1\}$.

It was shown by Cao and Chen \cite{Cao-Chen} that the tensor $B_{ijk}$ agrees with the Cotten tensor of $(M,g)$, up to a constant factor. In particular, we have $B_{ijk} = 0$ on the Bryant soliton. \\

\begin{proposition}
If $(M,g)$ is a steady gradient Ricci soliton, then we have 
\begin{align}
\label{B}
|B|^2 &= -(1-R) \, \Delta R - \frac{3}{4} \, |\nabla R|^2 \notag \\ 
&- \langle \nabla f,\nabla R \rangle - R^2 \, (1-R). 
\end{align} 
\end{proposition}

\textbf{Proof.} Using (\ref{gradient.of.R}) and (\ref{Laplacian.of.R}), we obtain  
\begin{align*} 
&\sum_{i,j,k} |B_{ijk}|^2 \\ 
&= \sum_{i,j,k} |\text{\rm Ric}_{ik} \, \partial_j f - \text{\rm Ric}_{ij} \, \partial_k f|^2 + \frac{1}{4} \, |\nabla R + 2R \, \nabla f|^2 \\
&- \sum_{i,j,k} (\text{\rm Ric}_{ik} \, \partial_j f - \text{\rm Ric}_{ij} \, \partial_k f) \, (\partial^j R + 2R \, \partial^j f) \, g^{ik} \\ 
&= \sum_{i,j,k} |\text{\rm Ric}_{ik} \, \partial_j f - \text{\rm Ric}_{ij} \, \partial_k f|^2 - \frac{1}{4} \, |\nabla R + 2R \, \nabla f|^2 \\ 
&= 2 \, |\text{\rm Ric}|^2 \, |\nabla f|^2 - 2 \, \sum_{i,j,k} \text{\rm Ric}_{ij} \, \partial_j f \, \text{\rm Ric}_{ik} \, \partial_k f - \frac{1}{4} \, |\nabla R + 2R \, \nabla f|^2 \\ 
&= 2 \, |\text{\rm Ric}|^2 \, |\nabla f|^2 - \frac{1}{2} \, |\nabla R|^2 - \frac{1}{4} \, |\nabla R + 2R \, \nabla f|^2 \\ 
&= -(\Delta R + \langle \nabla f,\nabla R \rangle) \, |\nabla f|^2 - \frac{1}{2} \, |\nabla R|^2 - \frac{1}{4} \, |\nabla R + 2R \, \nabla f|^2 \\ 
&= -|\nabla f|^2 \, \Delta R - |\nabla f|^2 \, \langle \nabla f,\nabla R \rangle - \frac{3}{4} \, |\nabla R|^2 \\ 
&- R \, \langle \nabla f,\nabla R \rangle - R^2 \, |\nabla f|^2. 
\end{align*} 
Using the identity $|\nabla f|^2 = 1-R$, we conclude that 
\begin{align*}
|B|^2 &= -(1-R) \, \Delta R - \frac{3}{4} \, |\nabla R|^2 \\ 
&- \langle \nabla f,\nabla R \rangle - R^2 \, (1-R), 
\end{align*} 
as claimed. \\

In the next step, we choose a smooth function $\psi: (0,1) \to \mathbb{R}$ such that $\nabla R + \psi(R) \, \nabla f = 0$ on the Bryant soliton. 

\begin{proposition}
If $(M,g)$ is a steady gradient Ricci soliton, then the vector field $X = \nabla R + \psi(R) \, \nabla f$ satisfies 
\begin{align} 
\label{general.formula}
&(1-R) \, \text{\rm div} \, X \notag \\ 
&= -|B|^2 - \frac{3}{4} \, \langle \nabla R - \psi(R) \, \nabla f,X \rangle \notag \\ 
&- \langle \nabla f,X \rangle + (1-R) \, \psi'(R) \, \langle \nabla f,X \rangle \\ 
&- \frac{3}{4} \, (1-R) \, \psi(R)^2 + (1-R) \, \psi(R) \notag \\ 
&- R^2 \, (1-R) + R \, (1-R) \, \psi(R) - (1-R)^2 \, \psi(R) \, \psi'(R). \notag
\end{align} 
\end{proposition}

\textbf{Proof.} 
Using (\ref{B}), we obtain 
\begin{align*} 
&(1-R) \, \text{\rm div} \, X \\ 
&= (1-R) \, \Delta R + (1-R) \, \psi(R) \, \Delta f + (1-R) \, \psi'(R) \, \langle \nabla f,\nabla R \rangle \\ 
&= -|B|^2 - \frac{3}{4} \, |\nabla R|^2 - \langle \nabla f,\nabla R \rangle - R^2 \, (1-R) \\ 
&+ (1-R) \, \psi(R) \, \Delta f + (1-R) \, \psi'(R) \, \langle \nabla f,\nabla R \rangle \\ 
&= -|B|^2 - \frac{3}{4} \, \psi(R)^2 \, |\nabla f|^2 - \frac{3}{4} \, \langle \nabla R - \psi(R) \, \nabla f,X \rangle \\ 
&- \langle \nabla f,X \rangle + \psi(R) \, |\nabla f|^2 \\ 
&- R^2 \, (1-R) + (1-R) \, \psi(R) \, \Delta f \\ 
&+ (1-R) \, \psi'(R) \, \langle \nabla f,X \rangle - (1-R) \, \psi(R) \, \psi'(R) \, |\nabla f|^2 \\ 
&= -|B|^2 - \frac{3}{4} \, (1-R) \, \psi(R)^2 - \frac{3}{4} \, \langle \nabla R - \psi(R) \, \nabla f,X \rangle \\ 
&- \langle \nabla f,X \rangle + (1-R) \, \psi(R) \\ 
&- R^2 \, (1-R) + R \, (1-R) \, \psi(R) \\ 
&+ (1-R) \, \psi'(R) \, \langle \nabla f,X \rangle - (1-R)^2 \, \psi(R) \, \psi'(R). 
\end{align*} 
This proves the assertion. \\

\begin{corollary} 
The function $\psi$ satisfies the differential equation 
\begin{equation} 
\label{ode}
0 = -\frac{3}{4} \, \psi(s)^2 + \psi(s) - s^2 + s \, \psi(s) - (1-s) \, \psi(s) \, \psi'(s)
\end{equation} 
for all $s \in (0,1)$.
\end{corollary}

\textbf{Proof.} 
The identity (\ref{general.formula}) holds for any steady gradient Ricci soliton. In particular, it holds for the Bryant soliton. On the other hand, we have $B = 0$ and $X = 0$ on the Bryant soliton. From this the assertion follows. \\

\begin{proposition} 
Assume that $\psi$ is chosen so that $\nabla R + \psi(R) \, \nabla f = 0$ on the Bryant soliton. Then 
\begin{align} 
\label{simplified.formula}
(1-R) \, \text{\rm div} \, X &= -|B|^2 - \frac{3}{4} \, \langle \nabla R - \psi(R) \, \nabla f,X \rangle \notag \\ 
&- \langle \nabla f,X \rangle + (1-R) \, \psi'(R) \, \langle \nabla f,X \rangle. 
\end{align} 
\end{proposition}

\textbf{Proof.} 
This follows immediately from (\ref{general.formula}) and (\ref{ode}). \\

In the next step, we consider the function 
\[u(s) = \log \psi(s) + \int_{\frac{1}{2}}^s \Big ( \frac{3}{2(1-t)} - \frac{1}{(1-t) \psi(t)} \Big ) \, dt.\] 

\begin{proposition}
We have 
\begin{align} 
\label{final.formula}
&(1-R) \, e^{-u(R)} \, \text{\rm div}(e^{u(R)} \, X) = -|B|^2 - \frac{R \, (R - \psi(R))}{\psi(R)^2} \, |X|^2. 
\end{align}
\end{proposition} 

\textbf{Proof.} 
Using (\ref{simplified.formula}), we obtain 
\begin{align*} 
&(1-R) \, e^{-u(R)} \, \text{\rm div}(e^{u(R)} \, X) \\ 
&= (1-R) \, \text{\rm div} \, X + (1-R) \, u'(R) \, \langle \nabla R,X \rangle \\ 
&= -|B|^2 - \frac{3}{4} \, \langle \nabla R - \psi(R) \, \nabla f,X \rangle \\ 
&- \langle \nabla f,X \rangle + (1-R) \, \psi'(R) \, \langle \nabla f,X \rangle \\ 
&+ \frac{3}{2} \, \langle \nabla R,X \rangle - \frac{1}{\psi(R)} \, \langle \nabla R,X \rangle + (1-R) \, \frac{\psi'(R)}{\psi(R)} \, \langle \nabla R,X \rangle \\ 
&= -|B|^2 + \frac{3}{4} \, \langle \nabla R + \psi(R) \, \nabla f,X \rangle - \frac{1}{\psi(R)} \, \langle \nabla R + \psi(R) \, \nabla f,X \rangle \\ 
&+ (1-R) \, \frac{\psi'(R)}{\psi(R)} \, \langle \nabla R + \psi(R) \, \nabla f,X \rangle \\ 
&= -|B|^2 + \frac{3}{4} \, |X|^2 - \frac{1}{\psi(R)} \, |X|^2 + (1-R) \, \frac{\psi'(R)}{\psi(R)} \, |X|^2. 
\end{align*} 
Using (\ref{ode}), we obtain 
\begin{equation} 
\label{ode.2}
-\frac{s \, (s - \psi(s))}{\psi(s)^2} = \frac{3}{4} - \frac{1}{\psi(s)} + (1-s) \, \frac{\psi'(s)}{\psi(s)}. 
\end{equation} 
Putting these facts together, the assertion follows. 

\section{Proof of the Theorem \ref{main.theorem}}

\label{main}

\begin{lemma}
\label{asymptotics.of.psi}
For $s \to 1$, we have $\psi(s) = \frac{2}{3} + O(\sqrt{1-s})$. 
\end{lemma} 

\textbf{Proof.} 
On the Bryant soliton, we have 
\[0 = \partial_i R + \psi(R) \, \partial_i f = -2 \sum_{i,j} \text{\rm Ric}_{ij} \, \partial_j f + \psi(R) \, \partial_i f.\] 
Therefore, the vector $\nabla f$ is an eigenvector of the Ricci tensor with eigenvalue $\frac{\psi(R)}{2}$. On the other hand, we have $\text{\rm Ric}_{ij} = \frac{1}{3} \, g_{ij} + O(|x|)$ near the origin. This implies $\frac{\psi(R(x))}{2} = \frac{1}{3} + O(|x|)$ near the origin. From this, the assertion follows easily. \\

\begin{lemma}
\label{asymptotics.of.u}
The limit $\lim_{s \to 1} u(s)$ exists.
\end{lemma}

\textbf{Proof.} 
It follows from Lemma \ref{asymptotics.of.psi} that 
\[\frac{1}{1-s} \, \Big ( \frac{3}{2} - \frac{1}{\psi(s)} \Big ) = O \Big ( \frac{1}{\sqrt{1-s}} \Big )\] 
for $s$ near $1$. Consequently, the limit $\lim_{s \to 1} (u(s) - \log \psi(s))$ exists. From this, the assertion follows. \\

\begin{proposition}
\label{inequality}
We have $\psi(s) < s$ for all $s \in (0,1)$.
\end{proposition}

\textbf{Proof.} 
Suppose that the assertion is false. Let 
\[s_0 = \sup \{s \in (0,1): \psi(s) \geq s\}.\] 
Since $\lim_{s \to 1} \psi(s) = \frac{2}{3}$, we conclude that $s_0 \in (0,1)$. Moreover, we have $\psi(s_0) = s_0$ and $\psi'(s_0) \leq 1$. Using (\ref{ode.2}), we obtain 
\begin{align*} 
0 &= -\frac{s_0 \, (s_0 - \psi(s_0))}{\psi(s_0)^2} \\ 
&= \frac{3}{4} - \frac{1}{\psi(s_0)} + (1-s_0) \, \frac{\psi'(s_0)}{\psi(s_0)} \\ 
&\leq \frac{3}{4} - \frac{1}{s_0} + \frac{1-s_0}{s_0} \\ 
&= -\frac{1}{4}. 
\end{align*} 
This is a contradiction. \\

\begin{proposition}
\label{div.theorem}
Let $\Omega$ be a bounded domain in $M$ with smooth boundary. Moreover, suppose that $R < 1$ at each point on $\partial \Omega$. Then 
\[\int_{\Omega \cap \{R<1\}} \frac{e^{u(R)}}{1-R} \, |B|^2 \leq -\int_{\partial \Omega} e^{u(R)} \, \langle X,\nu \rangle.\] 
\end{proposition}

\textbf{Proof.} 
Let us fix a smooth cut-off function $\chi: [0,\infty) \to [0,1]$ such that $\chi(s) = 0$ for $s \leq 1$ and $\chi(s) = 1$ for $s \geq 2$. Moreover, let $\chi_\varepsilon(s) = \chi \big ( \frac{s}{\varepsilon} \big )$. It follows from Proposition \ref{inequality} and (\ref{final.formula}) that 
\[(1-R) \, e^{-u(R)} \, \text{\rm div}(e^{u(R)} \, X) \leq -|B|^2.\] 
Using the divergence theorem, we obtain 
\begin{align} 
\label{a}
\int_{\partial \Omega} \chi_\varepsilon(1-R) \, e^{u(R)} \, \langle X,\nu \rangle 
&= \int_{\Omega \cap \{R<1\}} \text{\rm div}(\chi_\varepsilon(1-R) \, e^{u(R)} \, X) \notag \\ 
&= \int_{\Omega \cap \{R<1\}} \chi_\varepsilon(1-R) \, \text{\rm div}(e^{u(R)} \, X) \notag \\ 
&- \int_{\Omega \cap \{R<1\}} \chi_\varepsilon'(1-R) \ e^{u(R)} \, \langle X,\nabla R \rangle \\ 
&\leq -\int_{\Omega \cap \{R<1\}} \chi_\varepsilon(1-R) \, \frac{e^{u(R)}}{1-R} \, |B|^2 \notag \\ 
&- \int_{\Omega \cap \{R<1\}} \chi_\varepsilon'(1-R) \ e^{u(R)} \, \langle X,\nabla R \rangle. \notag
\end{align} 
We claim that 
\begin{equation}
\label{b}
\int_{\Omega \cap \{R<1\}} \chi_\varepsilon'(1-R) \ e^{u(R)} \, \langle X,\nabla R \rangle \to 0
\end{equation} 
as $\varepsilon \to 0$. Indeed, on the set $\Omega \cap \{\varepsilon \leq 1-R \leq 2\varepsilon\}$, we have 
\begin{align*} 
\chi_\varepsilon'(1-R) \, e^{u(R)} \, |\langle X,\nabla R \rangle| 
&\leq C_1 \, \varepsilon^{-1} \, |X| \, |\nabla R| \\ 
&\leq C_2 \, \varepsilon^{-1} \, (|\nabla R| + |\nabla f|) \, |\nabla R| \\ 
&\leq C_3 \, \varepsilon^{-1} \, |\nabla f|^2 \\ 
&= C_3 \, \varepsilon^{-1} \, (1-R) \\ 
&\leq 2C_3. 
\end{align*}
Here, $C_1$, $C_2$, and $C_3$ are positive constants which may depend on $\Omega$, but not $\varepsilon$. This implies 
\[\bigg | \int_{\Omega \cap \{R<1\}} \chi_\varepsilon'(1-R) \ e^{u(R)} \, \langle X,\nabla R \rangle \bigg | \leq 2C_3 \, \text{\rm vol}(\{\varepsilon \leq 1-R \leq 2\varepsilon\}),\] 
and the right hand side converges to $0$ as $\varepsilon \to 0$. This proves (\ref{b}). Combining (\ref{a}) and (\ref{b}), we conclude that 
\[\int_{\partial \Omega} e^{u(R)} \, \langle X,\nu \rangle \leq -\int_{\Omega \cap \{R<1\}} \frac{e^{u(R)}}{1-R} \, |B|^2,\] 
as claimed. \\

We now complete the proof of Theorem \ref{main.theorem}. By assumption, we can find an exhaustion of $M$ by bounded domains $\Omega_l$ such that 
\[\lim_{l \to \infty} \int_{\partial \Omega_l} e^{u(R)} \, \langle X,\nu \rangle = 0.\] 
Using Proposition \ref{div.theorem}, we obtain 
\[\int_{\Omega_l \cap \{R<1\}} \frac{e^{u(R)}}{1-R} \, |B|^2 \leq -\int_{\partial \Omega_l} e^{u(R)} \, \langle X,\nu \rangle.\] 
Passing to the limit as $l \to \infty$ gives 
\[\int_{\{R<1\}} \frac{e^{u(R)}}{1-R} \, |B|^2 = 0.\] 
Therefore, the tensor $B$ vanishes on the set $\{R < 1\}$. On the other hand, it is easy to see that the set $\{R<1\}$ is dense. Therefore, the tensor $B$ vanishes identically. It now follows from work of Cao and Chen \cite{Cao-Chen} that $(M,g)$ is rotationally symmetric.

\end{document}